%% file: LearningGBcomplexity.tex
\documentclass[12pt,noinfoline]{imsart}

\usepackage{amsmath,amssymb,amsthm,amscd,amsfonts,graphicx}%
\usepackage{fancyhdr, geometry}
\usepackage{mathtools,tikz}
\usepackage{url, hyperref, color, verbatim}
\usepackage[normalem]{ulem}
\usepackage{wrapfig, caption, subcaption,enumerate}
 \usepackage{natbib}
\theoremstyle{plain} 
\newtheorem{thm}{Theorem}[section]

\theoremstyle{definition}

\newtheorem{problem}[thm]{Problem}

\usepackage{tikz}
\usepackage{stmaryrd}

\def\jdlqed{\vbox{\hrule \hbox{\vrule\hbox to
5pt{\vbox to 6pt{\vfil}\hfil}\vrule}\hrule}}

\newcommand{\sonja}[1]{{\it\textcolor{magenta}{ [#1]}}}

\graphicspath{{./fig/}}

\begin{document}

\title{Learning a performance metric   of Buchberger's algorithm}
\runtitle{Learning Gr\"obner bases complexity}
\author{Jelena Mojsilovi\'c, Dylan Peifer, Sonja Petrovi\'c\thanks{JM is a PhD student at Purdue University.  SP is with Illinois Tech's Applied Math Department, and partially supported by the Simons Foundation Collaboration Grant for Mathematicians \texttt{854770}.  At the time of the first version of this manuscript, JM was an Applied Math undergraduate researcher at Illinois Tech and DP was a Mathematics PhD student at Cornell. }} 
\address{jmojsilo@purdue.edu, 
djp282@cornell.edu, sonja.petrovic@iit.edu}

\begin{abstract}
What can be (machine) learned about the performance of Buchberger's algorithm?

Given a system of polynomials, Buchberger's algorithm computes a Gr\"obner basis of the ideal these polynomials generate
using an iterative procedure based on multivariate long division. The runtime of each step of the algorithm is typically
dominated by a series of polynomial additions, and the total number of these additions is a hardware-independent
performance metric that is often used to evaluate and optimize various implementation choices. In this work we attempt
to predict, using just the starting input, the number of polynomial additions that take place during one run of
Buchberger's algorithm. Good predictions are useful for quickly estimating difficulty and understanding what features
make a Gr\"obner basis computation hard. Our features and methods could also be used for value models in the reinforcement
learning approach to optimize Buchberger's algorithm introduced in the second author's thesis. 

We show that a multiple linear regression model built from a set of easy-to-compute ideal generator statistics can
predict the number of polynomial additions somewhat well, better than an uninformed model, and better than regression
models built on some intuitive commutative algebra invariants that are more difficult to compute. We also train a simple
recursive neural network that outperforms these linear models. Our work serves as a proof of concept, demonstrating that
predicting the number of polynomial additions in Buchberger's algorithm is a feasible problem from the point of view of
machine learning.
\end{abstract}

\maketitle
\section{Introduction}

\emph{Does a given system of equations have a solution?} is a simple-to-state question in computational mathematics that has, perhaps surprisingly, proven generally difficult to answer. The question is ubiquitous in many areas of science, engineering, and mathematical applications, but mathematicians have ready-made scalable solutions to the problem only in the case of systems of linear equations: the linear polynomials are reduced to a basis from which one can deduce the answers. On the other hand, systems of non-linear multivariate polynomials have generated an amazing breadth of theoretical research in commutative algebra, algebraic geometry, and the computational branches of the field, including numerical.
Gr\"obner bases are algebraic structures that generalize the reduced row-echelon form for linear systems to ideals of multivariate polynomials.
 \cite{StWhatIsGB} offers a high-level overview of Gr\"obner bases and the textbook \cite{clo} discusses various applications as well. Of the many applications we single out two that have generated tremendous interest in recent decades: discrete optimization \cite{GBandOptim,rekha} and statistics \cite{GBandStats, St}. 

Many computer algebra systems offer generic algorithms for computing Gr\"obner bases applicable to all kinds of input ideals.  In the 1960s, Buchberger developed a groundbreaking algorithm \citep{Buchberger1965Translated} to compute a Gr\"obner basis of {\it any} ideal, a  problem that is $NP$-hard in general. As it applies to any polynomial system, Buchberger's algorithm has a doubly exponential  runtime in the number of variables \citep{dube}. In the decades that followed, several specialized algorithms have been used to improve runtime: \citet{BePa08,BePa08a,clo}. 
These algorithms form the cornerstone of the field of symbolic computational nonlinear algebra. 

  Improvements on  Buchberger's algorithm, such as Faug\`{e}re's  famous F5 algorithm developed in \cite{faugereetal, faugere2014sparse}, leverage the fact that the computation is a generalization of Gaussian elimination. As such, these methods construct nontrivial organizational techniques, such as cleverly organizing monomials into large matrices, to judiciously perform Buchberger's key step: reduction of S-polynomials. 
Namely, the algorithm grows a given generating set by adding nonzero remainders of S-polynomials upon division by the current generating set;  an S-polynomial is an element of the ideal created from a pair of given polynomials.  
The correctness of Buchberger's algorithm does not depend on the order in which such pairs, called S-pairs, are processed: one can simply create all pairs, store them in a queue in arbitrary order, and process them linearly.  On the other hand,  the above mentioned algorithms indicate -- and by now this is part of   the computational nonlinear algebra folklore --  one can improve the runtime by reorganizing the polynomials in the generating set so as to process the S-pairs in different order. Standard strategies select pairs based on some minimality criterion, such as the degree of the least common multiple of the two leading terms or its position in the monomial order.

In recent years,  symbolic computation and nonlinear algebra have been enriched with randomness. Randomization can be used to improve algorithm performance, as demonstrated, for example, in \cite{AlgebraicViolators, primality, spielmanteng1, BSKW}. 
Similarly,  machine learning can be used to predict various ingredients in  symbolic algebraic computations;  \cite{MLMathStructures} offers an overview of recent research on machine-learning mathematical structures; see also \cite{DeepLearnForSymbolicMath}. 
Machine learning (ML) offers a less-explored avenue for potential dramatic speed-ups in symbolic computation.  \cite{treeSearchFactorization} used tree search to choose the variable ordering for Horner factorization.  \cite{10.1007/978-3-319-08434-3_8} used support vector machine to pick a variable ordering for CAD. A review on why machine learning has particular challenges for symbolic computation is offered in \cite{10.1007/978-3-030-52200-1_29}, while \cite{10.1007/978-3-319-96418-8_20} gives an  overview of ML in mathematical software. 
In contrast to using ML to predict answers directly, \cite{LearnForDiscrete} uses ML to, among other things,  improve a step in the Simplex method. 
Within the paradigm of using learning to improve algorithms that give the exact answer, \cite{DylanMike-LearningBuch} uses machine learning to discover new S-pair selection
strategies in Buchberger's algorithm which outperform state-of-the-art human-designed heuristics by 20\% to 40\%.
Their main contribution was to express S-pair selection in Buchberger's algorithm as a reinforcement learning problem,
i.e., a game where the player or agent selects S-pairs and is rewarded for minimizing the overall computational cost of
the algorithm. Their measure of computational cost was the number of polynomial additions performed, which is a
hardware-independent number that indicates how hard the basis was to compute. A key part of many reinforcement learning
techniques is a value model which learns to predict future reward. For the setup of \cite{DylanMike-LearningBuch},
which did not use value models, such a model would predict the number of future polynomial additions before the
Gr\"obner basis computation is complete. The goal of this manuscript is to learn a version of this value function in
the supervised learning setting. 

\begin{problem}[General]\label{prob:general}
	For a given ideal $I$, how many polynomial additions are performed during one run of Buchberger's algorithm to compute a Gr\"obner basis of $I$? 
\end{problem}

Of course the problem as generally stated is not tractable - one needs to restrict to a type of ideal and a particular monomial ordering. Here we focus on the default ordering in {\tt Macaulay2} \citep{M2}, namely, graded-reverse-lexicographic.  
We consider two  families of ideals: binomial and toric, and generate random samples from each family.   
These types of ideals are ubiquitous in applications from integer programming \cite{gbfast} to statistics \cite{DS98}; a comprehensive reference is \cite{St}.  
When considering samples of ideals, it is important to note that `random data' is not the same as `generic data' in the commutative-algebra sense of the phrase; to this end, Section~\ref{sec:data} contains a discussion on avoidance of generic behavior in some detail, as well as the distribution of various features of the data sets we generate. 

The case for studying random binomial ideals has been well argued  in commutative algebra and has been also summarized in the second author's recent work \cite{DylanMike-LearningBuch} from the point of view of reinforcement learning: 
binomial ideals embody all the richness and complexities of computations with ideals; some of the hardest polynomial problems are binomial; and they can be generated randomly to avoid generic, or uninteresting, behavior, meaning avoiding zero-dimensional ideals\footnote{An ideal is said to be zero-dimensional if the solution of the generating polynomial system is a set of isolated points,  rather than a positive-dimensional variety. One expects the zero-dimensional case to be easier from the perspective of computing Gr\"obner bases.}.
Section~\ref{sec:data}  further illustrates that binomial ideals accurately capture much of the Gr\"obner basis problem: there is a
large variance in difficulty within distributions,  difficulty increases as expected when we increase the number of variables, and mostly as expected when we increase the number of generators of the ideal.  
We therefore define two instances of the general problem above.

\begin{problem}[Binomial]\label{prob:specific.bin}
Study Problem~\ref{prob:general} for binomial ideals in $3$ variables generated by 
	\begin{enumerate}[a)] 
	\item $4$  binomials of degree up to $20$, or 
	\item $10$ binomials of degree up to $20$. 
	\end{enumerate} 
\end{problem}

Toric ideals are a special class of binomial ideals  with a rich combinatorial structure; when defined by an integer design matrix $A$, the ideal is denoted by $I_A$. We postpone the formal definition until Section~\ref{sec:toric} and state the specific problem instance here:  
\begin{problem}[Toric]\label{prob:specific.toric}
	Study Problem~\ref{prob:general} for  toric ideals $I_A$ in $8$ variables obtained from non-negative integer matrices $A\subset \mathbb Z^{D\times 8}$ with:
	\begin{enumerate}[a)] 
	\item $D=2$ rows and integer entries up to $5$;
	\item $D=4$ rows and integer entries up to $5$;
	\item $D=6$ rows and integer entries up to $5$;
	\item $D=6$ rows and integer entries up to $10$.
	\end{enumerate} 	
\end{problem}

Note that in the toric problem, one does not know   the number of generators of  $I_A$ a priori.
We would be amiss not to mention that there exists fast software for efficiently computing both generators and Gr\"obner bases of $I_A$ from the input $A$: {\tt 4ti2} \citep{4ti2}.  
While one may not wish to use the standard implementation of Buchberger in these examples in practice, we have nevertheless found it   very interesting to see how well or how poorly the learning algorithm predicts the number of polynomial additions, as well as how distributions of various ideal features differ from random binomial ideals. 

\section{Ideal distributions} 
In commutative algebra, there are various approaches one can take for generating random samples of ideals. 
A set of polynomials can be generated in an ad hoc way from a given space, by taking, say, random monomials and adding them with random coefficients. Checking  distributions of features of the resulting ideals can help ensure the sample represents a varied enough set of polynomials for learning. While this way is heuristic, it offers a fast approach for obtaining a large amount of training samples. 

Another formal approach is to construct probabilistic models for polynomial ideals that can be used  to generate training samples with specific properties on demand. 
Such formal models require also understanding the induced distributions of various system invariants, and  are introduced in the third author's work \cite{rmi} on random monomial ideals.
  
Here, we generate datasets for the two problems above by considering the following models for  random binomial and toric ideals. 

\subsection{Random binomial ideals}





Considering an ideal generated by some polynomials, there are three basic parameters that are closely related to
computational difficulty: 
\begin{itemize}
\item $n$, the number of variables, 
\item $d$, the maximum degree of a monomial in the polynomials, and 
\item $s$, the number of
polynomial generators of the ideal.
\end{itemize}
 These parameters naturally map to the parameters of a random binomial ideal
model. First, choose $n$ as the fixed number of variables. Next, select $d$ as the maximal degree of any term in the
generating set. Finally, choose $s$ as the number of sampled binomials taken as generators of the ideal. We sample $s$
binomials by sampling $s$ pairs of two distinct monomials from the set of monomials with degree less than or equal to
$d$ in $n$ variables. There are two ways we considered sampling these monomials. The first, \emph{uniform}, samples two
distinct monomials uniformly at random from all monomials with degree less than or equal to $d$. The second,
\emph{weighted}, selects the degree of each monomial uniformly at random from 1 to $d$, then selects each monomial
uniformly at random among monomials of the chosen degree (in other words, weighting the monomials to ensure an equal
chance of each degree). The difference between these distributions is that weighted tends to produce more binomials of
low total degree while uniform samples a mix with mostly high degree. Both distributions assign non-zero coefficients
uniformly at random. We denote these distributions with the format ``$n$-$d$-$s$-(uniform/weighted)" to specify our
distribution on $s$-tuples of binomials of degree $\leq d$ in $n$ variables.

\subsection{Random toric ideals}\label{sec:toric}  
\input{toricideals.tex}

\subsection{The data}\label{sec:data}
\input{data}

\section{Linear regression models}\label{sec:regression}
\input{regression.tex}

\section{Learning via recursive neural networks (RNNs)}\label{sec:learning} 

Machine learning can be thought of as  a set of techniques used to infer  an underlying rule from a dataset. From a mathematical perspective, we think of such rules as different types of function fitting. The function can, for example, describe a relationship between variables, as in the regression setting. Neural networks are computational models within machine learning that originally drew inspiration from scientific models of biological neural networks. Each neural network defines a family of functions from which we can construct a rule for a given dataset. There are several parameters for each model that are tuned using the training data; details of this are omitted here, but 
\cite[\S 3.3.2]{DylanPhD} offers an overview of neural networks and further general references to machine learning. 

%


The main goal of this paper is to compare an out-of-the-box machine learning method with the naive human `guess' as well as the regression models constructed above. We have found  that the machine learning method does outperform regression, which is a proof of concept that properties like the ones we study in Problem~\ref{prob:specific.bin} are something that \emph{can be learned}. To illustrate, 
we first select one of the distributions of binomial ideals: the 3-20-10-weighted distribution. Recall that the dataset consists of a sample of size $1,000,000$, with $100,000$ samples set aside as a test set.

\begin{figure}[!h]
  \centering
  \includegraphics[width=\textwidth]{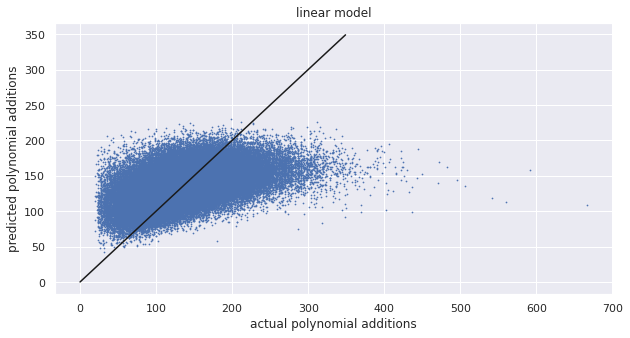}
  \includegraphics[width=\textwidth]{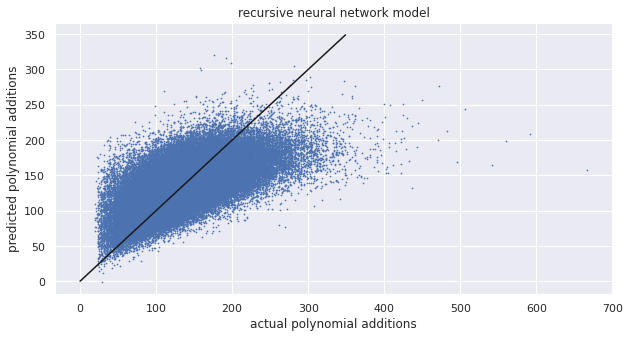}
  \caption{Predicted versus actual polynomial additions on a holdout test set of 100,000 examples in
    3-20-10-weighted. The dark line shows the location of perfect predictions.}
  \label{fig:eval-supervised}
\end{figure}

We will consider three basic models of increasing complexity. 
The first \emph{uninformed} model simply computes the mean value of polynomial additions on the training set and guesses that value for every input. This is equivalent to taking a simple `human' approach to prediction:  without any a priori knowledge of these distributions, simply guess the mean. 
Next, a \emph{linear} model
is trained on several hand-designed features (see previous section). 
In particular, we compute the maximum, minimum, mean, and standard
deviation of the degrees of the leading terms in the generating set, and also compute the number of pure power lead
terms (i.e., monomials that are a power of a single variable). These features are quick to compute and relate to the
algorithm difficulty.  
Finally, the most complex model we consider is a \emph{recursive neural network (RNN)} with gated recurrent unit (GRU) cells
containing $128$ hidden units followed by a single dense layer, for a total of $52,353$ parameters. The choice of the network topology (that is, the number of hidden units and a number of layers) is something that is generally adjusted ad hoc in any application; we started with some standard choices based on the work in \cite{DylanPhD}.   The network is trained
on a matrix composed of the exponent vectors of both terms on each generator, which is a $10 \times 6$ matrix of
integers. This allows the RNN to learn its own features from the input (i.e., the network is not given the max, min, mean, and standard deviation we used above in the linear models, but could conceivably compute them). Final agent performance is listed in
Table~\ref{tab:eval-supervised}, and the relationship between predicted and real values on the test set is shown in
Figure~\ref{fig:eval-supervised}.

\begin{table}[h]
  \centering
  \begin{tabular}{r|ccc}
    \hline
    model & Mean Squared Error & Mean Absolute Error & $R^2$ \\
    \hline
    uninformed & 2555.58 & 39.46 & 0.0000 \\
    linear & 1946.40 & 33.71 & 0.2384 \\
    RNN & 1499.37 & 29.25 & 0.4133 \\
    \hline
  \end{tabular}
  \caption{Trained model performance on a holdout test set of 100,000 examples in 3-20-10-weighted. The uninformed model
    simply guesses the mean on the training set, which is 135.44.}
  \label{tab:eval-supervised}
\end{table}

Both the linear and RNN model struggle on the extreme values of actual polynomial additions, as seen in the wide spread
of the plots in Figure~\ref{fig:eval-supervised}. They do show improvements over the uninformed model, and the RNN does
better than the linear model. For example, the mean absolute error shows that the RNN model averages a residual of 29.25
versus the uninformed model's 39.46 when predicting polynomial additions, of which an average value is
135.44. Unfortunately, the RNN model trained on 3-20-10-weighted has negative $R^2$ when evaluated on 3-20-4-weighted
and 3-20-4-uniform distributions, indicating that it is performing worse than guessing the mean on these related
distributions.

We also trained RNNs on the other distributions.
Table~\ref{table.rnn.R2} shows that the these learning models actually made some significant improvements; for example,  compare the values on diagonal of this table to those of Tables~\ref{table.TORregression.R2} and~\ref{table.Binregression.R2}.  
Unlike some of the regression models, RNNs trained on this data set did not generalize as well to other distributions, and due to the input format restriction, we were not able to evaluate RNN models trained on binomial data on toric data and vice versa. 

\begin{table}

\begin{tabular}{ c|c|c|c|c} 
Train$\backslash$Test data & 3-20-10-weighted & 3-20-10-uniform & 3-20-4-weighted & 3-20-4-uniform\\
\hline \hline 
 3-20-10-weighted & 0.41  & 0.06 & -4.58 & -6.36\\
 3-20-10-uniform & 0.23 & 0.13 & -5.06 & -6.62\\
 3-20-4-weighted & 0.06 & -0.63 & 0.37 & 0.19 \\
 3-20-4-uniform & 0.03  & -0.36 & 0.29 & 0.22\\ 
 \hline
\end{tabular} 
\medskip 

\begin{tabular}{ c|c|c|c|c}
Train$\backslash$Test data &   $\mathcal T(2,0,5,8)$  &  $\mathcal T(4,0,5,8)$   &  $\mathcal T(6,0,5,8)$ & $\mathcal T(6,0,10,8)$ \\
\hline \hline 
  $\mathcal T(2,0,5,8)$ & 0.39 & -0.04 & -15.43 & -0.32 \\ 
 $\mathcal T(4,0,5,8)$   & -2.55 & 0.15 & -78.20 & -3.03 \\ 
  $\mathcal T(6,0,5,8)$   & -2.63 & -0.00 & 0.65 & 0.04 \\ 
  $\mathcal T(6,0,10,8)$  & -3.31 & -0.07 & 0.05 & 0.52\\ 
  \hline
\end{tabular} 

  \caption{Summary of neural network predictions: $R^2$ for various training and testing data sets. Analogous to Tables~\ref{table.TORregression.R2} and~\ref{table.Binregression.R2}, but for the recursive neural network model, rather than the multiple linear regression model.
  }\label{table.rnn.R2}
\end{table}
\medskip 

\section{Concluding remarks} 

We set out to demonstrate that a performance metric of Buchberger's algorithm is something that \emph{can be predicted} from data, and that it is \emph{machine learnable}.  Prediction of this value  is intimately tied to reformulation of Buchberger's algorithm from the point of view of reinforcement learning. 

We have evaluated the performance of various linear regression models and a recursive neural network model (RNN) on random samples drawn from several distributions of toric and binomial ideals. We created these models and selected their parameter values to capture the expected behavior of Gr\"obner basis computations in practice. 
Linear regression models attained a two-fold goal: to create an approach that adequately models polynomial additions while remaining computationally efficient, and to use this model as a benchmark for RNN.

Several lines of future work can be derived directly from our explorations. 
Namely, there is a  need to continue studying other distributions of random ideals, defined so that parameters can be adjusted to achieve desired properties, similarly to what was done for monomial ideals in  \cite{rmi} or \cite{SerkanMFR}. 
Then, machine learning methods can be further fine-tuned to see how well they can perform on Problems~\ref{prob:general}, \ref{prob:specific.bin}, and \ref{prob:specific.toric}. This will likely require a computing infrastructure beyond what was available during our project. 
Finally, while binomial ideal data used in \cite{DylanMike-LearningBuch} has been shared in standard machine learning venues, we believe there is value in 
 collecting the data generated for various random polynomial explorations in a single location. To this end, the first and third author have collaborated with a team to implement  the starting steps of a Gr\"obner basis database \cite{GBdatabase} (see also \cite{JelenaURJ} for more information). 

The data we used can be downloaded from Zenodo at \url{10.5281/zenodo.6599502}, along with the code for learning and regression. The latter is also housed at \url{https://github.com/Sondzus/LearningGBvaluemodel}.

\bibliography{randomized-ideals,MLandAlg,SparkRandomizer}

\end{document}

%% file: toricideals.tex
A toric ideal is a prime binomial ideal defined as follows. 
Let 
\( A = \{a_1 \dots a_n\}\subseteq\mathbb Z^D\setminus\{0\} \)
be an integer matrix of rank $d\leq D$ with columns $\{a_1,\dots,a_n\}$.  
Denote by $t^{a_i}$ the monomial $t_1^{a_{1,i}}\cdots t_D^{a_{D,i}}$ in $k[t_1,\dots,t_D]$. 
Let $\varphi$ be the monomial map defined by $A$:
\begin{align*}
	\varphi: k[x_1,\dots,x_n] &\to k[t_1,\dots,t_D]  \\
	 x_i &\mapsto t^{a_i}. 
\end{align*} 
The toric ideal $I_A$ of the matrix $A$ is  the kernel  
\[I_A := \ker\varphi = \{ x^u-x^v : u-v\in\ker A\}.\]

Thus the natural way to generate toric ideals randomly is to generate the monomial map $\varphi$ 
randomly and then compute its kernel. 
To generate a random monomial map, we generate a random set of monomials using the  Erd\"os-Ren\'yi-type model for monomial ideals from \cite{rmi}. To allow for negative exponents (which can be worked around by homogenizing the matrix $A$ instead), we extend the random monomial ideal model to allow for Laurent monomials.   
The toric  model thus has 4 parameters: 
\begin{itemize}
\item $D$, the number of target variables, equiv. the number of rows of $A$;
\item $L$, the  bound on the negative total degree of each $\varphi(x_i)$: 
		$\sum_{\{j: a_{ij}<0\}} |a_{ij}| \leq L$; 
\item $U$, the  bound on the positive total degree of each $\varphi(x_i)$: 
		$\sum_{\{j: a_{ji}>0\}} a_{ji} \leq U$; 
\item $n$, the number of monomials, equiv. number of source variables, columns of $A$. 
\end{itemize} 
(The expert reader will notice that the matrices can be homogenized and the lower bound  set to $L=0$.)

Let us denote by $\mathcal T(D,L,U,n)$ the resulting distribution of toric ideals with these parameters. 
To compute the kernel $I_A$ as a starting basis for Buchberger's algorithm, we use \cite{FourTiTwo}, a {\tt Macaulay2} \citep{M2} interface for {\tt 4ti2} \citep{4ti2}, as it is the fastest software available to compute toric ideals. 

%% file: data.tex
We generated data sets from the following model specifications: 

\begin{center}
\begin{tabular}{ c|c|c} 
 model & type of ideals & sample size \\
\hline \hline 
 3-20-10-weighted& binomial  & 1,000,000 \\ 
 3-20-10-uniform& binomial  & 1,000,000  \\ 
 3-20-4-weighted& binomial  & 1,000,000\\ 
 3-20-4-uniform & binomial  & 1,000,000  \\ 
  $\mathcal T(2,0,5,8)$ & toric  & 429,093  \\ 
 $\mathcal T(4,0,5,8)$  & toric  & 314,688  \\ 
  $\mathcal T(6,0,5,8)$ & toric  &  325,927 \\ 
  $\mathcal T(6,0,10,8)$  & toric  & 151,532   \\ 
  \hline
\end{tabular} 
\end{center}
\medskip Note that toric ideal samples are smaller; this is because randomly generating binomials directly is very fast,
while obtaining a generating set of $I_A$ from $A$ is a difficult computation even for {\tt 4ti2}, so we just generated
as many samples as was computationally feasible.

\paragraph{Avoiding generic behavior.}

\begin{figure}[t]
  \centering
  \includegraphics[width=\textwidth]{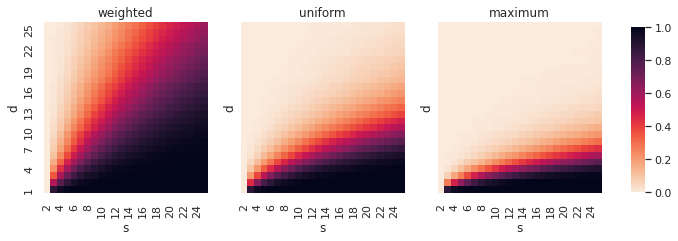}
  \caption{Proportion of zero-dimensional ideals in samples of 10,000 ideals each,  from \cite{DylanPhD}, for $n = 3$, $1 \le d \le 25$,
    $2 \le s \le 25$ and three model families (two of which we use here: weighted and uniform). Darker colors indicate more zero-dimensional ideals.} 
  \label{fig:dimension}
\end{figure}

\cite{DylanPhD} illustrates the avoidance of generic behavior in several random binomial ideal samples: 
Table~\ref{tab:dimension}  shows that all the six distributions of binomial ideals considered therein produce small proportions of zero-dimensional ideals.  Figure~\ref{fig:dimension} provides more in-depth information: increasing the number of generators increases the number of zero-dimensional ideals, and increasing the degree decreases the number of zero-dimensional ideals. 
\begin{table}[t]
  \centering
  \begin{tabular}
  {r|cc|cc}
    \hline
    & \multicolumn{2}{c|}{3-20-4} & \multicolumn{2}{|c}{3-20-10} \\
    $\dim$ & weighted & uniform 
    	& weighted & uniform  \\
    \hline
    0 & 188 & 4 
     & 2142 & 86 \\
    1 & 6263 & 2894 
	    & 7667 & 8198 \\
    2 & 3549 & 7102
     & 191 & 1716 \\
    \hline
  \end{tabular}
  \caption{Dimension of the binomial ideals in samples of 10,000 over  
  four 3-variable distributions from \cite{DylanPhD}.  } 
  \label{tab:dimension}
\end{table}

Additional illustrations in \cite{DylanPhD}, namely Figures 4.2 and 4.4 therein, further demonstrate binomial ideals as sufficient for reinforcement learning purposes: they show sizes of and maximum degrees of a generator in reduced Gr\"obner bases in samples of $10,000$ binomial ideals over six different $3$-variable distributions. 

\begin{figure}
\begin{centering}  \bf Distribution of invariants for binomial data from 3-20-10-uniform \end{centering} 
\smallskip

\begin{subfigure}{.3\linewidth}
\includegraphics[scale=0.4]{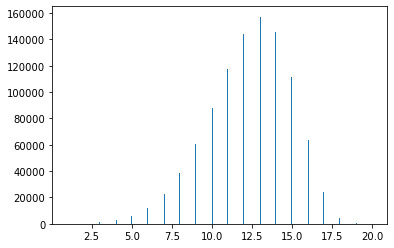}
\caption{Minimum generator degree}
\end{subfigure} 
\qquad   
\begin{subfigure}{.3\linewidth}
\includegraphics[scale=0.4]{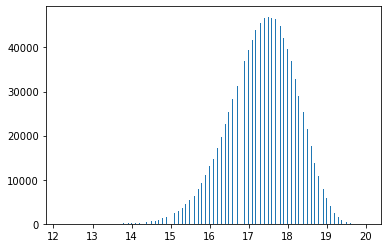}
\caption{Mean generator degree}
\end{subfigure} 
\qquad   
\begin{subfigure}{.3\linewidth}
\includegraphics[scale=0.4]{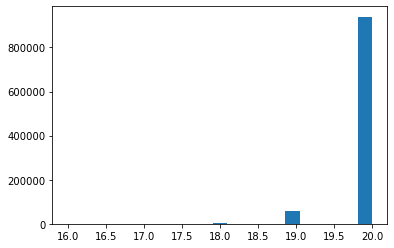}
\caption{Maximum  degree}
\end{subfigure} 

\begin{subfigure}{.3\linewidth}
\includegraphics[scale=0.4]{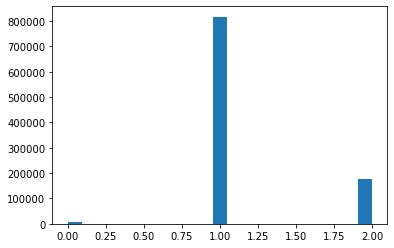}
\caption{Krull dimension}
\end{subfigure} 
\qquad   
\begin{subfigure}{.3\linewidth}
\includegraphics[scale=0.4]{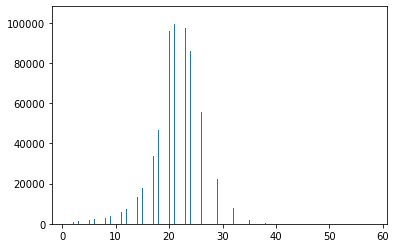}
\caption{Regularity}
\end{subfigure} 
\qquad   
\begin{subfigure}{.3\linewidth}
\includegraphics[scale=0.4]{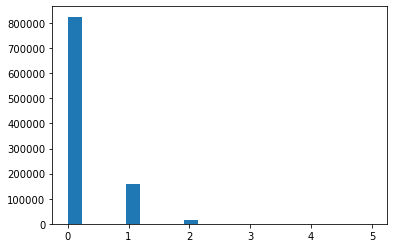}
\caption{Number of pure powers}
\end{subfigure} 

\caption{Histograms of sampling distributions of some features/invariants from the 3-20-10-uniform distribution of binomial ideals.
}\label{figure.binomial.histograms}
\end{figure}

\begin{figure}
\begin{centering}  \bf Distribution of invariants for toric data from  $\mathcal T(6,0,5,8)$\end{centering} 
\smallskip

\begin{subfigure}{.3\linewidth}
\includegraphics[scale=0.4]{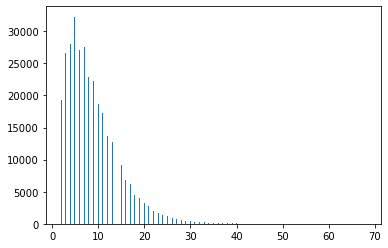}
\caption{Minimum generator degree}
\end{subfigure} 
\qquad   
\begin{subfigure}{.3\linewidth}
\includegraphics[scale=0.4]{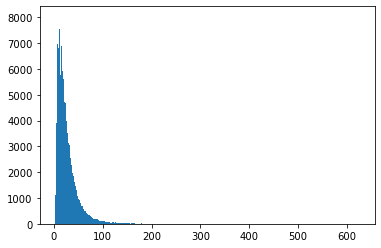}
\caption{Mean generator degree}
\end{subfigure} 
\qquad   
\begin{subfigure}{.3\linewidth}
\includegraphics[scale=0.4]{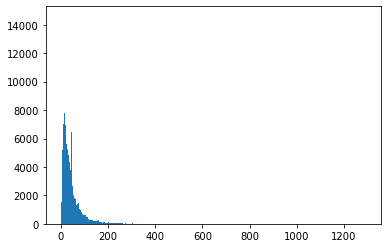}
\caption{Maximum  degree}
\end{subfigure} 

\begin{subfigure}{.3\linewidth}
\includegraphics[scale=0.4]{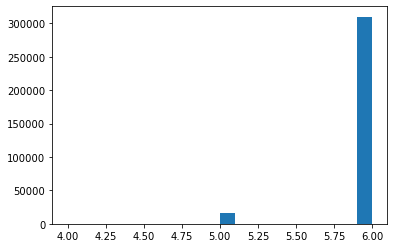}
\caption{Krull dimension }
\end{subfigure} 
\qquad   
\begin{subfigure}{.3\linewidth}
\includegraphics[scale=0.4]{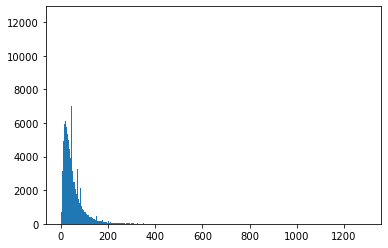}
\caption{Regularity}
\end{subfigure} 
\qquad   
\begin{subfigure}{.3\linewidth}
\includegraphics[scale=0.4]{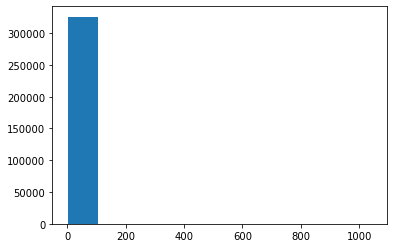}
\caption{Number of ideal generators}
\end{subfigure} 

\caption{Histograms of sampling distributions of some features/invariants from the $\mathcal T(6,0,5,8)$ distribution of toric ideals. Note the long right tails in most of the distributions.
}\label{figure.toric.histograms}
\end{figure}

\paragraph{How diverse are the data sets?} 
Figure~\ref{figure.binomial.histograms} illustrates sampling distributions of some basic algebraic invariants of samples of random binomial ideals: generator degree statistics (minimum, mean, maximum), Krull dimension and regularity of the coordinate ring, and the number of  pure powers appearing in the leading term of a generator.  For the non-expert reader, Krull dimension is the algebraic counterpart to the dimension of the solution set, and regularity   encodes the complexity of relations (or syzygies) on the generating polynomials.  
Figure~\ref{figure.toric.histograms} illustrates these invariants for random toric ideals, with the number of pure powers replaced by  the number of generators of $I_A$, which varies because the only the matrix $A$ is input to the sampling algorithm. 
 Admittedly, not all of these are reasonable invariants to be used as features in learning. For example, computing regularity is very expensive and as such it does not make sense to use it to predict Gr\"obner  bases complexity. On the other hand, regularity is a folklore `feature' which commutative algebraists use as an intuitive measure of  complexity of a Gr\"obner basis of an ideal. Perhaps surprisingly, we will discover that degree statistics for toric ideal generators outperform regularity in predicting the number of polynomial additions in a Gr\"obner basis computation. 

\begin{figure}
  \centering
  \includegraphics[width=.45\textwidth]{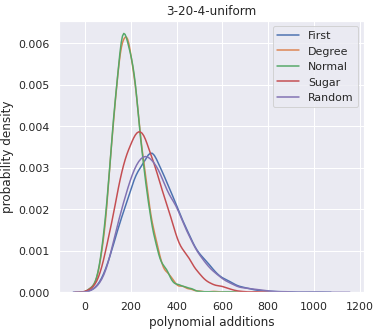}
  \includegraphics[width=.45\textwidth]{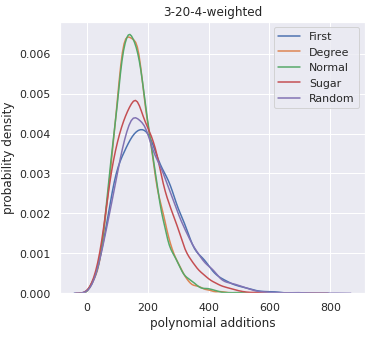}
  \caption{Kernel density estimations of  
  number of polynomial additions in 3-20-4 distributions over 10,000 samples. 
The figures show various S-pair selection strategies. \cite{DylanPhD}}
  \label{fig:kdes}
\end{figure}

\paragraph{How does the Buchberger algorithm perform on these data?}      
Examining just mean performance can be misleading, as there is high variance in difficulty within many of the
distributions; in particular, the distributions of the number of polynomial additions have large variance and long right tail \cite{DylanPhD}. Performance of various S-pair selection strategies (First, Degree, Normal, and Sugar) on random binomial ideals is illustrated in  Table~\ref{tab:additions}  and 
Figure~\ref{fig:kdes}, which indicates how better strategies both improve the mean and decrease the standard deviation. 
We do not discuss the specific selection strategies in detail here as they are not our main focus. 

\begin{table}[ht]
  \centering
  \begin{tabular}{r|cccc}
    \hline
    $n$ & First & Degree & Normal & Sugar \\
    \hline
    2 & 36.49[7.28] & 32.19[5.61] & 31.90[5.40] & 32.48[6.14] \\
    3 & 53.13[17.84] & 42.32[12.97] & 42.65[13.26] & 44.32[15.00] \\
    4 & 86.43[40.39] & 64.29[28.60] & 66.40[29.96] & 70.22[32.90] \\
    5 & 152.06[86.66] & 108.45[58.72] & 117.46[64.96] & 120.84[68.35] \\
    6 & 278.74[174.82] & 199.42[120.17] & 221.89[133.54] & 222.70[142.48] \\
    7 & 529.38[359.92] & 381.23[241.35] & 435.90[274.43] & 432.96[298.43] \\
    8 & 1042.88[784.55] & 767.37[516.61] & 890.31[588.71] & 874.39[654.15] \\
    \hline
  \end{tabular}
  \caption{Mean number of polynomial additions for different selection strategies on the same samples of 10,000
    ideals. Distributions are $n$-5-10-weighted. Table entries show mean[stddev]. \cite{DylanPhD}}
  \label{tab:additions}
\end{table}


\begin{figure}[!htbp]
  \centering
  \begin{subfigure}[b]{.80\textwidth}
    \centering
    \includegraphics[width=\textwidth]{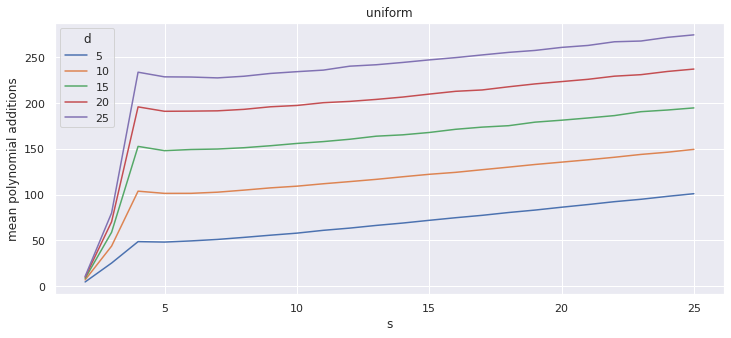}
  \end{subfigure}
  \begin{subfigure}[b]{.80\textwidth}
    \centering
    \includegraphics[width=\textwidth]{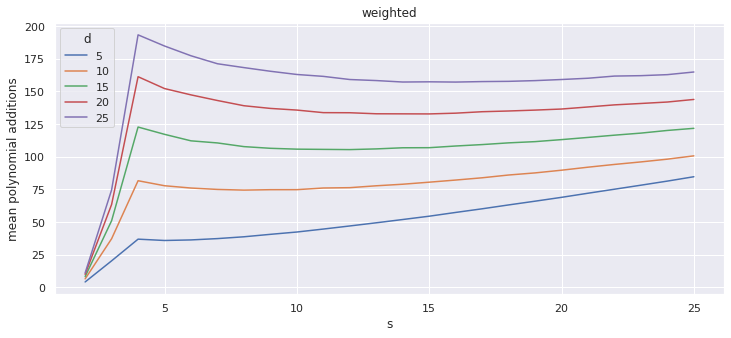}
  \end{subfigure}
  \caption{Average number of polynomial additions following Degree selection strategy in $n = 3$ as a function of $s$, for varying $d$ (see legend). Each degree and generator
    point is the mean over 10,000 samples. The distributions on the left are 3-d-s-uniform and on the right 3-d-s-weighted.}
  \label{fig:ds-additions}
\end{figure}

A few remarks are in order, extracted from results in \cite{DylanPhD}. The high variance observed in the distribution of the number of polynomial additions is a challenge for training reinforcement  learning models, yet it is exactly the behavior one desires the random ideal models to capture. Reflecting on Table~\ref{tab:additions}, one sees difficulty increasing rapidly with the number of variables, which is expected by the doubly-exponential worst-case run time of the algorithm.  The dependence on the number of generators $s$ is more subtle. Figure~\ref{fig:ds-additions} shows the previously noted spike in difficulty at $4$ generators and a slow increase after that. 
%

\begin{figure}[h]
  \centering
  \includegraphics[width=.5\textwidth]{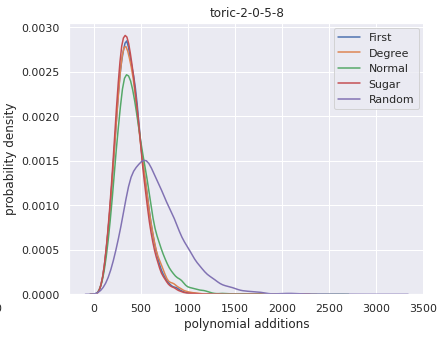}
  \caption{Kernel density estimations of the number of polynomial  additions in the 
   toric ideal distribution  $\mathcal T(2,0,5,8)$ over 10,000    samples.}
  \label{fig:toric-dist}
\end{figure}

Toric ideals generated from the probabilistic models we defined above can have a much wider range of degrees, generators, and difficulties than the binomial ideals, but carefully adjusting the parameters can lead to distributions with similar properties to our binomial model. Some examples are given in Figure~\ref{fig:toric-dist} and  Table~\ref{tab:toric-additions}.

\begin{table}[!htbp]
  \centering
  \begin{tabular}{r|cccc}
    \hline
    parameters & First & Degree & Normal & Sugar \\
    \hline
    2-0-5-8 & 390.55[155.49] & 395.45[162.32] & 437.21[199.85] & 382.38[151.27] \\
    4-0-5-8 & 416.31[545.45] & 400.30[485.37] & 404.93[490.27] & 396.61[483.46] \\
    6-0-5-8 & 31.64[71.28] & 31.00[57.03] & 31.01[57.06] & 30.97[56.92] \\
    6-0-10-8 & 68.36[161.92] & 68.42[163.28] & 68.42[163.28] & 68.39[162.57] \\
    \hline
  \end{tabular}
  \caption{Mean number of polynomial additions for different selection strategies on the same samples of 10,000 toric
    ideals. Table entries show mean[stddev].}
  \label{tab:toric-additions}
\end{table}

%% file: regression.tex

It makes little sense to deploy complex technology  such as machine learning to solve  Problems~\ref{prob:specific.bin} and~\ref{prob:specific.toric} if the number of polynomial additions  in Buchberger's algorithm can be predicted in easier ways. To that end, in this section we demonstrate the relative success of using simple and multiple linear regression for this prediction using the  data sets  from Section~\ref{sec:data}. 

In statistical modeling, the quantity we seek to predict  is called the \emph{response variable}. 
The response variable can be computed from a fixed, but unknown, function of the \emph{predictor variables}.  Fitting a model to the data means estimating this function; in linear regression modeling, the assumption is that the function is linear. As one may expect, linearity is rarely the `true model' in practice in applied statistics, however, linear regression models offer  surprisingly rich information about the data nonetheless. 
Predictor variables are selected by the modeler; the next sections discuss  types of predictors natural to use for Problem~\ref{prob:general}. 
Once the linear model is fit, one obtains coefficients for each of the predictors as well as a measure of how significant the predictor is. Roughly speaking, coefficients can be interpreted as a measure of unit change in the response per unit change in the predictor.  Significance  means that the given predictor is actually influential in the computation of the value of the response and  is measured using a $p$-value. Each predictor's  $p$-value for a given data is the probability of observing the given data, or more extreme, if the predictor is \emph{not} influential; thus low $p$-value shows evidence of influence. 

For both the binomial and toric ideals, we seek to use predictor variables that are easily computed, such as statistics on degrees of ideal generators, so that they can be readily used for prediction in practice. 
We also explore some of the  intuitive measures of  Gr\"obner complexity,  such as the Krull dimension, degree, and  Castelnuovo-Mumford regularity  of the coordinate ring of the ideal,  as they help in exploring how well regression can predict the number of polynomial additions.  
%
%

\subsection{Prediction using `easy' statistics} 
The following statistics of degrees of generators are natural choices to use  as predictors in linear regression: 
 the maximum, minimum, mean, and standard deviation  of generator degrees, and the  number of pure power leading terms of generators.  
In the toric ideal data set, by definition, the number of generators is not fixed, so this can also be  used as a predictor. Recall that the number of generators is fixed in the binomial data set. 

\tiny 
\begin{table}
\begin{tabular}{ c|c|c|c|c|c|c|c|c|c} 
 & MaxDeg & MinDeg & MeanDeg & StdDeg & PurePwrs & Deg & Dim & Reg & NumGens \\
\hline \hline 
 3-20-10-w & -0.34  & 6.84 & 11.69 & -9.75 &  -1.39 & -0.03 & -4.11 & 0.19 & -- \\ 
 3-20-10-u &  -1.69  & -0.79 & 10.65 & 2.27 &  1.42 & -0.03 & -10.13 & 0.80 & -- \\ 
 3-20-4-w & 2.05 & 3.13 & 9.15 & -3.14 & -12.14 & -0.18 & -21.91 & 0.16 & -- \\ 
 3-20-4-u & 1.84  & 	-0.02 {\tt *}   & 6.21 & -0.52 & -2.82 & -0.18 & -30.79 & 1.65 & -- \\ 
 \hline
  $\mathcal T(2,0,5,8)$ & 8.85  & -30.39 & 147.44 & -71.75 & -- & 3.70 & -- & 30.34 &  6.08 \\ 
 $\mathcal T(4,0,5,8)$  & -15.81 & -66.45  & 76.58 & -10.92 & -- &	-1.02 & -263.08 & 13.49 & 15.98 \\ 
  $\mathcal T(6,0,5,8)$ & -0.88 &  -0.55 & -1.13 & 1.97 & -- & 0.06 & -52.76 & 0.44 & 8.57 \\ 
  $\mathcal T(6,0,10,8)$  & -0.16 & -0.10  & 0.11 & 0.17 & -- & -0.00 & -384.88 & 0.08 & 4.33 \\ 
  \hline
\end{tabular} 
  \caption{Summary of regression analyses: Coefficients of various predictors in multiple linear regression for each of the 4 random binomial data sets and 4 random toric ideal data sets. For the entry marked with a {\tt *} (MinDeg in the  3-20-4-u dataset), the $p$-value was $0.501$. All other predictors are significant with  $p$-value less than $10^{-3}$. 
Data names were shortened for space considerations:  "w" and "u" stand for weighted and uniform, respectively. 
}\label{table.multiplelinearregression.coef}
\end{table}
\medskip 
\normalsize 

Table~\ref{table.multiplelinearregression.coef} shows that as dimension increases,  the number of polynomial additions are greatly reduced for all datasets, except $\mathcal T(2,0,5,8)$, but we will see later in Tables~\ref{table.regression.MMMSDDegBin} and~\ref{table.regression.MMMSDDegTor}  that this is because the $R^2$ value is $0$.    Other promising candidates for regression predictors are the minimum, maximum, mean, and standard deviation of generator degrees for all of the data sets; number of pure powers in the binomial data; and number of generators in the toric data. 
One interesting trend, or rather lack thereof, to note is that predictors are not always positively or negatively correlated with the response; for example, maximum degree has  a negative correlation with the number of polynomial additions in 3-20-10-weighted, but a positive correlation in 3-20-4-weighted. 
All of the predictors are  statistically significant with the exception of minimum generator degree, as shown in the entry for 3-20-4-uniform; we thus re-train the linear regression model without it when performing model validation in Section~\ref{sec:validate}. 

Note that using output of the type  listed in Table~\ref{table.multiplelinearregression.coef} alone does not indicate which predictors are best at predicting polynomial additions; ultimately, one builds different models and compares their performance, as we illustrate next. 

\subsection{
Exploring intuitive predictors}  

Krull dimension, degree, and Castelnuovo-Mumford regularity are natural predictors a commutative algebraist might consider using.  While these are in general difficult to compute, except for the dimension of a toric ideal which can be obtained from the defining matrix, bounds on them can be used as indicators of Gr\"obner basis complexity. 
We explore how linear models based on these individual predictors perform in terms of correctly estimating the number of polynomial additions in Buchberger's algorithm. Table~\ref{table.regression.MMMSDDegBinTor} offers a summary of the comparison of these models with a 
  linear regression model built with less computationally expensive predictors: Minimum, Maximum, Mean, and Standard Deviation Degree, which we abbreviate as  MMMSDDeg. 

\tiny 
\begin{table}
\begin{tabular}{ c|c|c|c|c} 
  & 3-20-10-weighted & 3-20-10-uniform & 3-20-4-weighted & 3-20-4-uniform \\
\hline \hline 

 MMMSDDeg & 0.248  & 0.017 & 0.169 & 0.029  \\ 
 Regularity & 0.041 & 0.010 & 0.035 & 0.036  \\ 
 Dimension & 0.010 & 0.003 & 0.000 & 0.014  \\ 
 Degree & 0.005 & 0.000 & 0.006 & 0.005 \\ 
 \hline
\end{tabular} 

\smallskip

\begin{tabular}{c|c|c|c|c} 
  & $\mathcal T(2,0,5,8)$  &  $\mathcal T(4,0,5,8)$   &  $\mathcal T(6,0,5,8)$ & $\mathcal T(6,0,10,8)$ \\
\hline \hline 

 MMMSDDeg & 0.226 & 0.278 & 0.129 & 0.572 \\ 
 Regularity & 0.264 & 0.156 & 0.051 & 0.148 \\ 
 Dimension & 0.000 & 0.000 & 0.214 & 0.018 \\ 
 Degree & 0.094  & 0.059 & 0.020 & 0.003 \\ 
 \hline
\end{tabular} 
\smallskip
  \caption{Summary of fitting the number of polynomial additions using various linear regression models: $R^2$ statistics for the multiple linear regression model built with the generator degree statistics (minimum, maximum, mean, standard deviation -- MMMSDDeg), and for each simple linear regression model built with each of regularity, dimension, and degree. 
  }\label{table.regression.MMMSDDegBinTor}
\end{table}
\medskip 
\normalsize 

The various linear models are compared using the $R^2$ statistic, or  the coefficient of determination. This is a statistical measure of how close the data lie to the estimated regression line. The baseline model which always predicts the mean will have an $R^2$ value of $0$, and the model that exactly matches the observed values will have $R^2=1$. Models that do worse than then baseline prediction will have a negative $R^2$ value; this often happens when the linear regression is trained (or fitted) on a particular data set and tests or evaluated on a completely different one (e.g., data draws from a completely different distribution). Table~\ref{table.regression.MMMSDDegBinTor} compares the $R^2$ values of the linear regression model built using MMMSDDeg with that of  simple linear regression models built on each of regularity, dimension, and degree. Often times the MMMSDDeg has a higher $R^2$ than any of the `natural' predictors.
Therefore, one can construct reasonably good models using the generator degree statistics which are easy to compute. 
Perhaps  somewhat surprisingly, it even often outperforms regularity in predicting the number of polynomial additions in a Gr\"obner computation for random toric ideals (particularly, note the last column in the Table!).

\tiny 
\begin{table}
\begin{tabular}{c|c|c|c|c} 
  & $\mathcal T(2,0,5,8)$  &  $\mathcal T(4,0,5,8)$   &  $\mathcal T(6,0,5,8)$ & $\mathcal T(6,0,10,8)$ \\
\hline \hline 
 number of generators  & 0.036  & 0.413 & 0.487 & 0.855 \\ 
 \hline 
\end{tabular} 
\medskip

\begin{tabular}{c|c|c|c|c} 
  & 3-20-10-w  & 3-20-10-u & 3-20-4-w & 3-20-4-u \\ 
  \hline \hline
 number of pure powers & 0.022 & 0.001 & 0.023 & 0.001 \\ 
 \hline
\end{tabular} 
\smallskip

  \caption{The $R^2$ statistics summarizing modeling the number of polynomial additions using  simple linear regression on number of generators in toric ideals, and on number of pure powers in generator lead terms in binomial ideals. 
  }\label{table.regression.NGPP}
\end{table}
\medskip 
\normalsize 

We also note that in the toric dataset, the number of generators of the ideal 
  appears to be quite useful in predicting polynomial additions, thus  we use it in the regression model to evaluate against other toric ideals. 
  In the binomial dataset, the number of pure powers in the lead terms of the generators has low $R^2$, see Table~\ref{table.regression.NGPP} for a summary. Still, we chose to use this predictor in regression in order to evaluate against other binomial datasets, since it is easily computable and its coefficient in multiple linear regression was not low in comparison to other coefficients while also being statistically significant. 

\subsection{Model validation} \label{sec:validate} 
As is customary in statistics, the regression models are trained  using 90\% of the data within each training data set. In other words, the training subset of the data is used to estimate the coefficients on each predictor variable and constant in linear regression.  
%
%
%
Estimated coefficients are listed in Table~\ref{table.multiplelinearregression.coef}; as mentioned earlier, coefficients that are \emph{not} significant are marked with a {\tt *}. 

The remaining 10\% of each dataset is used to test each linear model's ability to predict the number of polynomial additions. The readers less familiar with statistics should note  that this is a methodologically correct way of validating the performance of the learning algorithm. 
It is not expected that the models perform well at all on distributions on which they were not trained; nevertheless, we tested them across the different distributions. What we found is that  regression models trained on binomial ideal data tend to be mildly generalizable to other binomial ideal distributions, but poorly so to toric ideal data in Table~\ref{table.Binregression.R2}. Likewise, models trained on toric ideal distributions do sometimes generalize to other toric distributions (for example, notice some of the relatively high $R^2$ values in the off-diagonal entries in the lower-right corner of Table~\ref{table.TORregression.R2}), but do not have a good ability to predict polynomial additions in binomial ideal distributions. 

Tables~\ref{table.TORregression.R2} and~\ref{table.Binregression.R2} summarize these findings and hint at the problem difficulty, but also show some promise of generalizability within a problem class to `nearby' distributions (those with similar parameter values).

\begin{table}
\begin{tabular}{ c|c|c|c|c} 
 Train$\backslash$Test data & 3-20-10-weighted & 3-20-10-uniform & 3-20-4-weighted & 3-20-4-uniform \\
\hline \hline 
 3-20-10-weighted & 0.24  & -0.02 & 0.001 & -0.02 \\ 
 3-20-10-uniform & -0.26 & 0.02 & 0.13 & 0.03 \\ 
 3-20-4-weighted & -0.09 & 0.003 & 0.17 & 0.02  \\ 
 3-20-4-uniform & -0.27  & 	0.02  & 0.14 & 0.03 \\ 
 \hline
  $\mathcal T(2,0,5,8)$ & -1482.65  & -1848.3 & -801.41 & -1172.01 \\ 
 $\mathcal T(4,0,5,8)$  & -62.73 & -31.14  & -8.04 & -4.32 \\ 
  $\mathcal T(6,0,5,8)$ & -2.00 &  -5.83 & -5.45 & -7.72 \\ 
  $\mathcal T(6,0,10,8)$  & -3.54 & -7.63  & -5.00 & -6.84 \\ 
  \hline
\end{tabular} 
  \caption{Summary of regression predictions: $R^2$ for various training and testing data sets. The rows in the table are indexed with training data, that is, the data used to fit the model and estimate the regression coefficients. The columns in the table are indexed by testing data, that is, the data used to evaluate how well model predicts the number of polynomial additions on that particular data set.  
  }\label{table.TORregression.R2}
\end{table}
\medskip 

\begin{table}
\begin{tabular}{ c|c|c|c|c} 
 Train$\backslash$Test data &  $\mathcal T(2,0,5,8)$  &  $\mathcal T(4,0,5,8)$   &  $\mathcal T(6,0,5,8)$ & $\mathcal T(6,0,10,8)$ \\
\hline \hline 
 3-20-10-weighted &  -4.54 & -0.12 & -23.18 & -193.91 \\ 
 3-20-10-uniform & -3.77 & -0.11 & -35.05 & -243.98 \\ 
 3-20-4-weighted & -4.64 & -0.09 & -47.06 & -423.88 \\ 
 3-20-4-uniform & -3.89 & -0.09 & -42.87 & -361.33 \\ 
 \hline
  $\mathcal T(2,0,5,8)$ & 0.25 & -0.34 & -9414.59 & -76130.1 \\ 
 $\mathcal T(4,0,5,8)$  & -4.78 & 0.81 & -889.82 & -8275.07 \\ 
  $\mathcal T(6,0,5,8)$ & -3.78 & 0.57 & 0.64 & -0.61 \\ 
  $\mathcal T(6,0,10,8)$ & -4.48 & 0.22 & 0.52 & 0.88\\ 
  \hline
\end{tabular} 
  \caption{Summary of regression predictions: $R^2$ for various training and testing data sets. The rows in the table are indexed with training data, that is, the data used to fit the model and estimate the regression coefficients. The columns in the table are indexed by testing data, that is, the data used to evaluate how well model predicts the number of polynomial additions on that particular data set.  
  }\label{table.Binregression.R2}
\end{table}
\medskip 

